\documentclass{amsart}[12pt]
\newcommand{\privatecomment}[1]{} 
\usepackage[applemac]{inputenc} 
\usepackage[all]{xy}
\usepackage{url}
\newcommand{\refequ}[1]{$(\ref{#1})$}

\theoremstyle{plain}
\newtheorem{thm}{Theorem}[section]

\newtheorem{lemma}[thm]{Lemma}
\newtheorem{prop}[thm]{Proposition}

\theoremstyle{definition}
\newtheorem{defin}[thm]{Definition}

\theoremstyle{remark}

\numberwithin{equation}{section}

\newcommand{\Or}{\mathcal{O}}
\newcommand{\hn}{n/2}
\newcommand{\BQ}{\mathbb Q}
\newcommand{\Bk}{\mathbf {k}}

\newcommand{\Apl}{\operatorname{A}_{\textrm{PL}}}
\newcommand{\coker}{\mathrm{coker}}
\newcommand{\im}{\mathrm{im}\,}

\newcommand{\Top}{\operatorname{Top}}
\newcommand{\CDGA}{\operatorname{CDGA}}
\newcommand{\Ho}{\operatorname{H}}
\newcommand{\Hoch}{\operatorname{HH}}

\begin{document}

\title[Poincar\'e duality and CDGA's]{Poincar\'e duality and commutative differential graded algebras}
\author{Pascal Lambrechts}
\author{Don Stanley}%
\address{P.L.: Universit\'e de Louvain, Institut Math\'ematique\\
2, chemin du Cyclotron\\B-1348 Louvain-la-Neuve, BELGIUM}
\email{lambrechts@math.ucl.ac.be}
\address{D.S.: University of Regina, Department of Mathematics\\
College West, Regina, CANADA}%
\email{stanley@math.uregina.ca}%
\thanks{P.L. is Chercheur Qualifi\'e au F.N.R.S.}%
\subjclass{55P62, 55M05, 57P10}
\keywords{Poincar\'e duality, Commutative differential graded algebra.}%
\date{}
\maketitle
\begin{abstract}
We prove that every commutative differential graded algebra
whose cohomology is a simply-connected Poincar\'e
duality algebra is quasi-isomorphic to one whose
underlying algebra is simply-connected and satisfies Poincar\'e
duality in the same dimension. This has applications in rational homotopy,
giving  Poincar\'e duality at the cochain level,
which is of interest in particular in the study of configuration spaces
and in string topology.
\end{abstract}
\vspace{1cm}

\begin{center}
\begin{minipage}{10cm}
\begin{center}
\textbf{Dualité de Poincaré et algèbres différentielles graduées commutatives.}\\
\end{center}
{\small\footnotesize\textsc{Résumé.} Nous démontrons que toute algèbre différentielle graduée commutative (ADGC) dont la cohomologie est une algèbre simplement connexe à 
dualité de Poincaré est quasi-isomorphe à une ADGC dont l'algèbre sous-jacente est à dualité de Poincaré dans la même dimension.
Ce résultat a des applications en théorie de l'homotopie rationnelle, permettant d'obtenir la dualité de Poincaré au niveau des cochaines,
entre autres dans l'étude des espaces de configurations et en topologie des cordes.}
\end{minipage}
\end{center}
\section{Introduction}
\label{section-introduction}

The first motivation for the main result of this paper comes from rational homotopy theory.
Recall that Sullivan \cite{Sul:inf} has constructed a contravariant functor
\[\Apl\colon\Top\to\CDGA_\BQ\]
from the category of topological space to the category of commutative differential graded algebras
over the field $\BQ$  (see Section \ref{sec-term} for the definition).
The main feature of $\Apl$ is that when $X$ is a simply-connected space with rational homology
of finite type, then the rational homotopy type of $X$ is completely encoded in any CDGA
 $(A,d)$ weakly equivalent to $\Apl(X)$. By \emph{weakly equivalent} we mean that
 $(A,d)$ and $\Apl(X)$ are connected by
 a zig-zag of CDGA morphisms inducing isomorphism in homology, or
 \emph{quasi-isomorphisms} for short,
 \[(A,d)\stackrel{\simeq}\leftarrow\dots\stackrel{\simeq}\rightarrow \Apl(X).\]
We then say that $(A,d)$ is a \emph{CDGA-model} of $X$
 (see \cite{FHT} for a complete exposition of this theory).
  We are particularly interested in the case when $X$
  is a simply-connected closed manifold of dimension $n$,
  since then $\Ho^*(A,d)$ is a simply-connected Poincar\'e duality algebra of dimension $n$
  (a graded algebra $A$ is said to be \emph{simply-connected} if $A^0$ is isomorphic to the ground field and
  $A^1=0$; see also Definition \ref{def-PDA}).

  Our main result is the following:

  \begin{thm}
\label{thm-main} Let $\Bk$ be a field of any characteristic  and let $(A,d)$
 be a CDGA over $\Bk$ such that $\Ho^*(A,d)$ is a
simply-connected Poincar\'e duality algebra in dimension $n$. Then there exists a
CDGA $(A',d')$ weakly equivalent to
$(A,d)$ and such that $A'$ is
a simply-connected algebra satisfying Poincar\'e duality in
dimension $n$.
\end{thm}
\par
Our theorem was conjectured by Steve Halperin over 20 years ago.
The $A'$ of the theorem is called a \emph{differential Poincar\'e duality algebra}
or \emph{Poincar\'e duality CDGA}
(Definition \ref{def-DPDA}). In particular any simply-connected closed manifold
 admits a Poincar\'e duality CDGA-model.

Notice that the theorem is even valid for a field of non-zero
characteristic. Also our proof is very constructive: Starting from
a finite-dimensional CDGA $(A,d)$, it shows how to compute
explicitly a weakly equivalent differential Poincar\'e duality
algebra $(A',d')$. We will also prove in the last section that
under some extra connectivity hypotheses, any two such weakly
equivalent differential Poincar\'e duality algebra can be
connected by a zig-zag of quasi-isomorphisms between differential
Poincar\'e duality algebras.

Aubry, Lemaire, and Halperin \cite{ALH} and Lambrechts \cite[p.158]{Lambrechts-cochainthick}
prove the main result of this paper in some special cases.
Also in \cite{Stasheff} Stasheff proves some chain level results about
Poincar\'e duality using Quillen models. An error in Stasheff's paper was corrected in
\cite{ALH}.

Before giving the idea of the proof of Theorem \ref{thm-main}, we
describe a few applications of this result.

\subsection{Applications}
There should be many applications of this result to constructions in rational homotopy theory
involving Poincar\'e duality spaces. We consider here two:
The first is to the study of configuration
spaces over a closed manifold, and the second to string topology.

Our first application is to the determination of the rational homotopy type of the configuration space
\[F(M,k):=\{(x_1,\dots,x_k)\in M^k:x_i\not=x_j\textrm{ for }i\not=j\}\]
of $k$ points in a closed manifold $M$ of dimension $n$. When $k=2$ and $M$ is $2$-connected, we showed in
\cite[Theorem 1.2]{LS-FM2} that if $A$ is a Poincar\'e duality CDGA-model of $M$
then a CDGA-model of $F(M,2)$ is given by
\begin{equation}\label{E:FA2}
A\otimes A/(\Delta)
\end{equation}
where $(\Delta)$ is the differential ideal in $A\otimes A$ generated by the so-called
\emph{diagonal class} $\Delta\in (A\otimes A)^n$.

For $k\geq2$ we have constructed in \cite{LS-dgmodFMk} an explicit CDGA
\[F(A,k)\]
generalizing \refequ{E:FA2} and which is an \emph{$A^{\otimes k}$-DGmodule model}
of $F(M,k)$.
Poincar\'e duality of the CDGA $A$ is an essential ingredient in the construction of
$F(A,k)$.
If $M$ is a smooth complex projective variety then we can use $H^*(M)$ as a model for
$M$, and in this case $F(H^*(M),k)$ is exactly the model of Kriz and
Fulton-Mac Pherson \cite{Kriz}\cite{FultonMacPherson}
for $F(M,k)$.
However we do not know in general if $F(A,k)$ is also a CDGA-model although it seems to be the natural
candidate.

\vspace{5mm}
A second application is to string topology, a new field created
by Chas and Sullivan \cite{ChasSullivan}.
They constructed a product, a bracket and a $\Delta$ operator
on the homology of the free loop space $LM=M^{S^1}$
of a closed simply-connected manifold $M$, that turned it into Gerstenhaber
algebra and even a BV algebra.
On the Hochschild cohomology $\Hoch^*(A,A)$ of a (differential graded)
algebra $A$, there are the classical cup product and Gerstenhaber bracket, and
Tradler \cite{Tradler-BV} showed that for $A=C^*(M)$ there is
also a $\Delta$ operator on Hochschild homology making it into a BV algebra.
Menichi \cite{Menichi} later reproved this result and showed that the
$\Delta$ can be taken to be the dual of the Connes boundary operator.
Recently Felix and Thomas \cite{FelixThomas} have shown that over the rationals
the Chas-Sullivan BV structure on
the homology of $LM$ is isomorphic to the BV structure on
$\Hoch^*(C^*(M),C^*(M))$. Their proof uses the main result of this paper.
Yang \cite{Yang} also uses our results to give explicit formulas for
the BV-algebra structure on Hochschild cohomology.

\subsection{Idea of proof}

The proof is completely constructive. We start with a CDGA $(A,d)$
and an orientation $\epsilon\colon A^n \rightarrow \Bk$
(Definition \ref{def-o}). We consider the pairing at the chain level
$$
\phi\colon A^k\otimes A^{n-k} \rightarrow \Bk, \ a\otimes b \mapsto \epsilon(ab)
$$
We may assume that $\phi$ induces a non degenerate bilinear form on cohomology
making $\Ho^*(A)$ into a Poincar\'e duality algebra.
The problem is that $\phi$ itself may be degenerate; there may be some orphan elements
(see Definition \ref{def-orient}) $a$ with $\epsilon(ab)=0$ for all $b$.
Quotienting out by the orphans $\Or$ we get a differential Poincar\'e duality algebra
$A/\Or$, and a map
$f\colon A\rightarrow A/\Or$ (Proposition \ref{prop-PDbarA}).
With this observation the heart of the proof begins.
\par
Now the problem is that $f$ might not be a quasi-isomorphism - this happens whenever
$\Ho^*(\Or)\not=0$. The solution is to add generators to $A$ to get a quasi-isomorphic algebra
$\hat A$ with better properties.
An important observation is that $\Ho^*(\Or)$ satisfies a kind of Poincar\'e duality so it is enough to
eliminate $\Ho^*(\Or)$ starting from about half of the dimension and working up from there.
In some sense we perform something akin to surgery by eliminating the
cohomology of the orphans in high dimensions
and having the lower dimensional cohomology naturally disappear at the same time.
In the middle dimension, the extra generators have the effect of turning orphans which represent
homology classes into non orphans. In higher dimensions some of the new generators become orphans
whose boundaries kill elements of $\Ho^*(\Or)$. In both cases the construction introduces
no new orphan homology between the middle dimension and the dimension where the
elements of $\Ho^*(\Or)$ are killed.
This together with the duality in $\Ho^*(\Or)$ is enough to get an inductive
proof of Theorem \ref{thm-main}.

\section{Some terminology}\label{sec-term}

Just for the record we introduce the terms CDGA, Poincar\'e duality algebra and
differential Poincar\'e duality algebra.

We fix once for all a ground field $\Bk$ of any characteristic. So tensor product, algebras, etc.,
will always be over that field.
A commutative differential graded algebra, or CDGA, $(A,d)$ is a non-negatively graded commutative algebra,
together with a differential $d$ of degree $+1$.
If an element $a\in A$ is in degree $n$, we write $|a|=n$. The set of elements of degree $n$ in $A$
is denoted $A^n$. Since $A$ is graded commutative we have the
formula $ab=(-1)^{|a||b|}ba$ and $a^2=0$ when $|a|$ is odd, including when $\Bk$ is
of characteristic $2$. Also $d$ satisfies the graded Leibnitz formula
$d(ab)=(da)b+(-1)^{|a|}adb$. CDGA over the rationals are of particular interest since
they are models of rational homotopy theory. For more details see \cite{FHT}.
\par
\noindent
{\bf Convention:}
All of the CDGA
we consider in this paper will be connected, in other words $A^0=\Bk$, and of finite type.

Note that every simply connected CW-complex of finite type admits such a
CDGA model of its rational homotopy type.

Poincar\'e duality is defined as follows:

\begin{defin}\label{def-PDA}
An {\em oriented Poincar\'e duality algebra} of dimension $n$ is a pair
$(A,\epsilon)$ such that
$A$ is a connected graded commutative algebra and $\epsilon\colon A^n\rightarrow \Bk$
is a linear map
such that the induced bilinear forms
$$
A^k\otimes A^{n-k} \rightarrow \Bk, \ a\otimes b \mapsto \epsilon(ab)
$$
are non-degenerate.
\end{defin}

The following definition comes from \cite{LS-FM2}:
\begin{defin}
\label{def-DPDA} An {\em oriented differential Poincar\'e duality
algebra} or {\em oriented Poincar\'e CDGA} is a triple $(A,d,\epsilon)$ such that
\begin{itemize}
\item[(i)] $(A,d)$ is a CDGA,
\item[(ii)] $(A,\epsilon)$ is an oriented Poincar\'e duality
algebra,
\item[(iii)] $\epsilon(dA)=0$
\end{itemize}
\end{defin}

An oriented differential Poincar\'e duality algebra is essentially a CDGA whose underlying algebra
satisfies Poincar\'e duality. The condition $\epsilon(dA)=0$ is equivalent to
$\Ho^*(A,d)$ being a Poincar\'e duality algebra in the same dimension
\cite[Proposition 4.8]{LS-FM2}.

For convenience we make the following:

\begin{defin}\label{def-o}
An \emph{orientation} of a CGDA $(A,d)$ is a linear map
$$\epsilon\colon A^n\to \Bk$$
such that $\epsilon(dA^{n-1})=0$ and there exists a cocycle
$\mu\in A^n\cap\ker d$ with  $\epsilon(\mu)=1$.
\end{defin}
Recall that $s^{-n}\Bk$ is the chain complex which is non-trivial only in degree $n$ where it is
$\Bk$.
Notice that the above definition is equivalent to the fact that
$\epsilon\colon(A,d)\to s^{-n}\Bk$ is a chain map that induces an
epimorphism $\Ho^n(\epsilon)\colon \Ho^n(A,d)\to\Bk$. We will use this alternative definition of
orientation interchangeably with the first without further comment.
The definition of differential Poincar\'e duality algebra can be thought of as a combination
of Definitions \ref{def-o} and \ref{def-PDA}.

If $V$ is a vector space and $v_1, \dots , v_l$ are elements of $V$, we let
$\langle v_1, \dots, v_l \rangle$ or $\langle \{ v_i \} \rangle$ denote the linear subspace spanned by
these elements.

\section{The set of orphans}
In this section we consider a fixed CDGA $(A,d)$ such that
$\Ho^*(A,d)$ is a connected Poincar\'e duality algebra in dimension
$n$.

The proof of our main theorem will be based on the study of orphans,
which is the main topic of this section.
\begin{defin}\label{def-orient}
If $\epsilon$ is an orientation on $(A,d)$ then the \emph{set of
orphans} of $(A,d,\epsilon)$ is the set
$$\Or:=\Or(A,d,\epsilon):=
\{a\in A |\forall b\in A, \epsilon(a\cdot b)=0\}.
$$
\end{defin}
\begin{prop}
The set of orphans $\Or$ is a differential ideal in $(A,d)$
\end{prop}
\begin{proof}
$\Or$ is clearly a vector space since $\Or=\cap_{b\in
A}\ker\left(\epsilon(b\cdot-)\right)$.

If $a\in\Or$ and $\xi\in A$ then for any $b\in A$ we have
$\epsilon((a\xi)b)=\epsilon(a(\xi b))=0$.
Therefore $\Or$ is an ideal.

If $a\in\Or$ then for any $b\in A$ we have, using the
fact that $\epsilon(dA)=0$,
$$\epsilon((da)b)=\pm\epsilon(d(ab))\pm\epsilon(a(db))=0.$$
Therefore $d\Or\subset\Or$.
\end{proof}
Clearly $\Or\subset\ker\epsilon$ since $\epsilon(\Or\cdot 1)=0$. Thus
the orientation $\epsilon\colon A\to s^{-n}\Bk$ extends to a chain map
$\bar\epsilon\colon \bar A:=A/\Or\to s^{-n}\Bk$ that also induces
an epimorphism in $\Ho^n$, and so $\bar \epsilon\colon \bar A \rightarrow s^{-n}\Bk$ is
itself an orientation.
\begin{prop}\label{prop-nada}
\label{prop-PDbarA} Let $(A,d)$ be a CDGA such that $\Ho^*(A,d)$ is
a Poincar\'e duality algebra in dimension $n$ and let
$\epsilon\colon A^n\to\Bk$ be an orientation. Assume that $A$ is
connected and of finite type. Let $\Or$ be the set of orphans of
$(A,d,\epsilon)$, let $(\bar A,\bar d):=(A,d)/\Or$, and let
$\bar\epsilon\colon\bar A\to s^{-n}\Bk$ be the induced
orientation.

Then $(\bar A,\bar d,\bar\epsilon)$ is an oriented differential
Poincar\'e duality algebra and $H(\bar A,\bar d)$ is a Poincar\'e
duality algebra in degree $n$.
\end{prop}
\begin{proof}
We know that $\bar\epsilon(d\bar A)=0$ since $\bar\epsilon$ is a chain map.
As in \cite[Definition 4.1]{LS-FM2} consider the bilinear form
$$\langle-,-\rangle\colon\bar A\otimes\bar A\to\Bk\,,\,\bar a
\otimes\bar b \mapsto\bar\epsilon(\bar a.\bar b)
$$
and the induced map
$$\theta\colon \bar A\to\hom(\bar A,\Bk)\,,\,\bar
a\mapsto\langle\bar a,-\rangle.$$ Let $\bar a=a\mod\Or\in \bar
A\setminus\{0\}$. Then $a\in A\setminus \Or$ and there exists
$b\in A$ such that $\epsilon(a.b)\not=0$. Set $\bar
b=b\mod\Or\in\bar A$. Then $\theta(\bar a)\not=0$ because
$\theta(\bar a)(\bar b) \not=0$. Thus $\theta$ is injective and
since $\bar A$ and $ \hom(\bar A,\Bk)$ have the same dimension
this implies that $\theta$ is an isomorphism and $(\bar A,\bar
d,\bar\epsilon)$ is a differential Poincar\'e duality algebra in
the sense of \ref{def-DPDA}. By \cite[Proposition 4.7]{LS-FM2}
$\Ho^*(\bar A,\bar d)$ is a Poincar\'e duality algebra in dimension
$n$.
\end{proof}

\begin{lemma}
\label{lemma-Okerd-dA} $\Or\cap\ker d\subset d(A)$.
\end{lemma}
\begin{proof}
Let $\alpha\in\ker d$ of degree $k$. If $\alpha\not\in\im d$ then
$[\alpha]\not=0$ in $\Ho^*(A,d)$ and by Poincar\'e duality there
exists $\beta\in A\cap \ker d$ of degree $n-k$ such that $[\alpha].[\beta]\not=0$.
Therefore $\epsilon(\alpha.\beta)\not=0$ and $\alpha\not\in\Or$.
\end{proof}
Consider the following short exact sequence
\begin{equation}
\label{equ-SES}
 \xymatrix{
0\ar[r]&\Or\ar@{^(->}[r]&A\ar[r]^-{\pi}&\bar A=A/\Or\ar[r]&0.
 }
 \end{equation}
 Notice that, in spite of Lemma \ref{lemma-Okerd-dA}, the differential
 ideal $\Or$ is in general not acyclic. When it is then the map
 $\pi$ is a quasi-isomorphism and Proposition \ref{prop-PDbarA}
 shows that $\bar A$ is the desired differential Poincar\'e duality
 model of $A$.
 The idea of the proof of our main theorem will be to modify $A$
 in order to turn the ideal of orphans into an acyclic ideal.
 Actually a Poincar\'e duality argument shows that its enough to
 get the acyclicity of $\Or$ in a range of degrees above half
 the dimension. In order to make this statement precise we introduce the following
 definition:
\begin{defin}
The set of orphans $\Or$ is said to be \emph{$k$-half-acyclic} if
$\Or^i\cap\ker d\subset d(\Or^{i-1})$ for $\hn+1\leq i\leq k$.
\end{defin}
In other words $\Or$ is $k$-half-acyclic iff $\Ho^i(\Or,d)=0$ for
$\hn+1\leq i\leq k$. Clearly this condition is empty for $k\leq
\hn$. Therefore an orphan set is always ($\hn$)-half-acyclic.
\begin{prop}
\label{prop-piqi} If $\Or$ is $(n+1)$-half-acyclic and $A$ is
connected and of finite type then $\pi\colon A\to \bar
A:=A/\Or$ is a quasi-isomorphism.
\end{prop}
\begin{proof}
By hypothesis $\Ho^*(A)$ is a Poincar\'e duality algebra in
dimension $n$ and by Proposition \ref{prop-PDbarA} the same is
true for $\Ho^*(\bar A)$. Moreover since $A^0=\Bk$, these cohomologies
are connected and $\pi^*=\Ho^*(\pi)$ sends the fundamental class of
$\Ho^n(A)$ to the fundamental class of $\Ho^n(\bar A)$. All of this
implies that $\pi^*$ is injective.

Thus the short exact sequence \refequ{equ-SES} gives us short
exact sequences
$$\xymatrix{
0\ar[r]& \Ho^i(A)\ar[r]^{\pi^*}&\Ho^i(\bar A)\ar[r]&(\coker
\pi^*)^i\cong \Ho^{i+1}(\Or)\ar[r]&0. }$$

By $(n+1)$-half-acyclicity, $\Ho^i(\Or)=0$ for $\hn+1\leq i\leq
n+1$. Also $(\coker \pi^*)^{>n}=0$. Thus $(\coker
\pi^*)^{\geq\hn}=0$. By Poincar\'e duality of $\Ho^*(A)$ and
$\Ho^*(\bar A)$ we deduce that $(\coker \pi^*)^{\leq\hn}=0$.
Therefore $\coker \pi^*=0$ and $\pi$ is a quasi-isomorphism.
\end{proof}
\section{A certain extension of a given oriented CDGA}
The aim of this section is,  given an integer $k\geq\hn+1$ and  an
oriented CDGA $(A,d,\epsilon)$, to construct a certain
quasi-isomorphic oriented CDGA $(\hat A,\hat d,\hat\epsilon)$. In
the next section we will prove that if the set $\Or$ of orphans of
$A$ is $(k-1)$-half-acyclic then the set $\hat\Or$ of orphans of
$\hat A$ is $k$-half-acyclic.

In this section we will always suppose that $(A,d)$ is a CDGA
equipped with a chain map $\epsilon\colon (A,d)\to s^{-n}\Bk$
satisfying the following hypotheses:\\
\begin{equation}\label{equ-2red}
\left\{
\begin{tabular}{l}
(i) $A$ is of finite type\\
(ii) $A^0\cong\Bk$, $A^1=0$, $A^2\subset \ker d$\\
(iii) $\Ho^*(A,d)$ is a Poincar\'e duality algebra in dimension
$n\geq 7$\\
(iv) $\epsilon\colon (A,d)\to s^{-n}\Bk$ is an orientation.
\end{tabular}
 \right.
 \end{equation}

We also suppose given a fixed integer $k\geq\hn+1$.

Next we start the construction of the oriented CDGA $(\hat
A,\hat d,\hat\epsilon)$. Set $l:=\dim(\Or^k\cap\ker
d)-\dim(d(\Or^{k-1}))$. Choose $l$ linearly independent elements
$\alpha_1,\ldots,\alpha_l\in\Or^k\cap \ker d$ such that
\begin{equation}\label{equ-Pascal}
\Or^k\cap \ker d=d(\Or^{k-1})\oplus\langle
\alpha_1,\ldots,\alpha_l\rangle.
\end{equation}
In a certain sense the $\alpha_i$'s are the obstruction to
$\Ho^k(\Or)$ being trivial. By Lemma \ref{lemma-Okerd-dA} there
exist $\gamma'_1,\ldots,\gamma'_l\in A^{k-1}$ such that
$d\gamma'_i=\alpha_i$.

Choose a  family $h_1,\ldots, h_m\in A\cap\ker d$ such that
$\{[h_i]\}$ is a homogeneous basis of $\Ho^*(A,d)$. Using the Poincar\'e
duality of $\Ho^*(A,d)$ there exists another family $\{h^*_i\}\subset
A\cap\ker d$ such that $\epsilon(h^*_j.h_i)=\delta_{ij}$, where
$\delta_{ij}$ is the Kronecker symbol. We set
$$
\gamma_i:=\gamma'_i-\sum_j \epsilon(\gamma'_i.h_j).h^*_j
$$
and
\begin{equation}\label{equ-Gamma}
\Gamma:=\langle\gamma_1,\ldots,\gamma_l\rangle\subset A^{k-1}.
\end{equation}
The two main properties of this family are the following:
\begin{lemma}
\label{lemma-G.kerd} $d(\gamma_i)=\alpha_i$ and
$\epsilon(\Gamma\cdot \ker d)=0$.
\end{lemma}
\begin{proof}
The first equation is obvious since $d\gamma'_i=\alpha_i$ and
$h^*_j$ are cocycles.

A direct computation shows that $\epsilon(\gamma_i\cdot h_j)=0$. On the
other hand using the facts that $\epsilon(\im d)=0$ and
$\alpha_i\in\Or$ we have that for $\xi\in A$,
$$
\epsilon(\gamma_i\cdot d\xi)=\pm\epsilon(d(\gamma_i\cdot \xi))\pm\epsilon(\alpha_i\cdot \xi)=0.
$$
Since $\ker d=\langle h_1,\ldots,h_m\rangle\oplus \im d$, the lemma
has been proven.
\end{proof}

Next using the above data we construct a relative Sullivan algebra $(\hat
A,\hat d)$ that is quasi-isomorphic to $(A,d)$ and with some new generators $c_i$ that bound
the $\alpha_i$. To define $(\hat A,\hat d)$ properly we distinguish two cases:
\par
\noindent
Case 1: when $char(\Bk)=0$ or $k$ is odd,
\begin{equation}\label{equ-Ahat}
(\hat A,\hat d):=(A\otimes\wedge(c_1,\ldots,c_l,w_1,\ldots,w_l);
\hat d(c_i)=\alpha_i,\hat d(w_i)=c_i-\gamma_i)
\end{equation}
Case 2: when $char(\Bk)$ is a prime $p$ and $k$ is even,
$$
(\hat A,\hat d):=(A\otimes\wedge(\{c_i,w_i,u_{i,j},v_{i,j}\}_{1\leq i\leq l,j\geq 1})
$$
with differential given by:
$$
\hat d(c_i)=\alpha_i,\hat d(w_i)=c_i-\gamma_i, \hat d(u_{i,1})=w_i^p,
\hat d u_{i,j}=v^p_{i,j-1},
$$
$$
\hat dv_{i,1}=(c_i-\gamma_i)w_i^{p-1},
\hat dv_{i,j}=v_{i,j-1}^{p-1}\hat dv_{i,j-1}
$$
Notice that $\deg(c_i)=k-1$, $\deg(w_i)=k-2$, $\deg(u_{i,1})=p(k-2)-1$ and
$\deg(v_{i,1})=p(k-2)$. All the other generators $u_{i,j}$ and $v_{i,j}$ have degree larger than $n$.
It will turn out that
only $c_i,w_i, u_{i,1}$ and $v_{i,1}$ will be relevant and the last two only when $k$ is small and
$p=2$.

\begin{lemma}
\label{lemma-jqi} The injection $j\colon(A,d)\to(\hat A,\hat d)$
is a quasi-isomorphism.
\end{lemma}
\begin{proof}
The lemma follows since the cofibre
$\hat A\otimes_A \Bk$ of $j$ is $\wedge(c_1,\ldots,c_l,w_1,\ldots,w_l)$ or
$\wedge(\{c_i,w_i,u_{i,j},v_{i,j}\}_{1\leq i\leq l,j\geq 1}$ which are acyclic.
\end{proof}

Our next step is to build a suitable orientation $\hat \epsilon$
on $\hat A$ that extends $\epsilon$. We construct this orientation
so that the $c_i$ are orphans (except when $k$ is about half the
dimension which requires a special treatment). This will prevent
the $\alpha_i$ from obstructing the set of orphans from having
trivial cohomology in  degree $k$.  In order to define this
orientation $\hat\epsilon$ we first need to define a suitable
complementary subspace of $\im d$ in $A$.

Next we choose a complement $Z$ of $\Or\cap d(A)$ in $\Or$.
\begin{lemma}
\label{lemma-ZG} $d(Z)=d(\Or)$, $Z\cap\Gamma=0$ and
$(Z\oplus\Gamma)\cap d(A)=0$.
\end{lemma}
\begin{proof}
The proof that $d(Z)=d(\Or)$ is straightforward.

Let $\gamma=\sum_ir_i\gamma_i\in Z\cap\Gamma$. Since
$Z\subset\Or$, $\alpha:=\sum r_i\alpha_i=d\gamma\in d(\Or)$.
By equation (\ref{equ-Pascal}) this implies that each $r_i=0$, hence $\gamma=0$ and
$Z\cap\Gamma=0$.

Let $z\in Z$ and $\gamma=\sum_ir_i\gamma_i\in \Gamma$. Suppose
that $z+\gamma\in\im d$. Then $d(z+\gamma)=0$, hence $\alpha:=\sum
r_i\alpha_i=d\gamma=-dz\in d(\Or)$. Again this implies that each
$r_i=0$ and $\gamma=0$. Therefore $z\in\im d$. By the definition of
$Z$ this implies that $z=0$.
\end{proof}

Choose a complement $U$ of $Z\oplus\Gamma\oplus d(A)$ in $A$. Set
$$T:=Z\oplus\Gamma\oplus U$$
which is a complement of $d(A)$ in $A$.

We are now ready to define our extension $\hat\epsilon$ on $\hat
A$. For $\xi\in A^+$ and $t\in T$ we set
\begin{equation}
\label{equ-epshat} \left\{
\begin{tabular}{l}
(i) $\hat\epsilon(\xi)=\epsilon(\xi)$\\
(ii) $\hat\epsilon(w_id(\xi))=(-1)^k\epsilon(\gamma_i\xi)$\\
(iii) $\hat\epsilon(c_ic_j)=-\epsilon(\gamma_i\gamma_j)$\\
(iv)
$\hat\epsilon(w_i)=\hat\epsilon(w_it)=\hat\epsilon(c_i)=\hat\epsilon(c_i\xi)=\hat\epsilon(c_ic_j\xi)=
\hat\epsilon(c_iw_j)=$\\
$\qquad\hat\epsilon(c_iw_j\xi)=\hat\epsilon(w_iw_j)=\hat\epsilon(w_iw_j\xi)=\hat\epsilon(u_{i,1})=
\hat\epsilon(w_iu_{j,1})=$\\
$\qquad\hat\epsilon(v_{i,1})=\hat\epsilon(u_{i,1}\xi)=\hat\epsilon(v_{i,1}\xi)=0$\\
(v) $\hat\epsilon(x)=0$ if $\deg(x)\not=n$.
\end{tabular} \right.
\end{equation}
\begin{lemma}

The formulas $(\ref{equ-epshat})$ define a unique linear map
$\hat\epsilon\colon\hat A\to s^{-n}\Bk$.
\end{lemma}
\begin{proof}
Let $x\in \hat A^n$.
Since $n\geq 7$, and $|w_i|, |c_i|\geq n/2-1$, for degree reasons
$x$ is of length at most $2$ in the $w_i$ and
$c_i$. Similarly $x$ is of length at most
one in $v_{i,1}$ and $u_{i,1}$. Moreover for $j>1$,
$|v_{i,j}|>|u_{i,j}|>n$, and $v_{i,1}w_j$, $v_{i,1}c_j$, and $u_{i,1}c_j$ all
have degree $>n$.
Also $A=T\oplus d(A)$. From these facts it follows that
$(\ref{equ-epshat})$ defines $\hat \epsilon$ on each monomial of $\hat A$. We can extend linearly to
all of $\hat A$.

Notice that $(\ref{equ-epshat})$(ii) is well defined since $\epsilon(\Gamma\cdot \ker d)=0$
by Lemma \ref{lemma-G.kerd}. Again using the fact that $A=T\oplus d(A)$, the well definedness
of $\hat\epsilon$ follows.

\end{proof}
\begin{lemma}\label{lemma-ori}
$\hat\epsilon\colon\hat A\to s^{-n}\Bk$ is an orientation.
\end{lemma}
\begin{proof}
We need to check that $\hat\epsilon(d(\hat A^{n-1}))=0$. Using the
fact that $\hat A^{\leq 1}=\hat A^0=\Bk$ and that $k\geq\hn+1$ we
get that every element of $\hat A^{n-1}$ is a linear combinations
of terms of the form $\xi$, $w_i\xi$, $c_i\xi$ for some $\xi\in A$
and possibly terms of the form $w_iw_j$, $w_ic_j$, $c_ic_j$,
$u_{i,1}\xi$ and $v_{i,1}$.
Using the definition \refequ{equ-epshat} of $\hat \epsilon$ we
compute:
\begin{itemize}
\item $\hat\epsilon(d\xi)=\epsilon(d\xi)=0$.
\item
$\hat\epsilon(d(w_i\xi))=
\hat\epsilon(c_i\xi)-\hat\epsilon(\gamma_i\xi)+(-1)^{\deg(w_i)}\hat\epsilon(w_id\xi)=0$
by formulas (iv) and (ii) of \refequ{equ-epshat}.
\item
$\hat\epsilon(d(c_i\xi))=\hat\epsilon(\alpha_i\xi)\pm\hat\epsilon(c_id\xi)=0$
because $\alpha_i\in\Or$.
\item
$\hat\epsilon(d(w_ic_j))=\hat\epsilon(c_ic_j)-\hat\epsilon(\gamma_ic_j)+(-1)^k\hat\epsilon(w_i\alpha_j)=
-\hat\epsilon(\gamma_i\gamma_j)+(-1)^k\hat\epsilon(w_i\alpha_j)
=-\hat\epsilon(\gamma_i\gamma_j)+\hat\epsilon(\gamma_i\gamma_j)=0$.
\item
$\hat\epsilon(d(c_ic_j))=\hat\epsilon(\alpha_ic_j)\pm\hat\epsilon(c_i\alpha_j)=0$.
\item
$\hat\epsilon(d(w_iw_j))=
\hat\epsilon(c_iw_j)-\hat\epsilon(\gamma_iw_j)\pm\hat\epsilon(w_ic_j)\pm\hat\epsilon(w_i\gamma_j)=0$
because $\gamma_i,\gamma_j\in T$.
\item
$\hat\epsilon(du_{i,1}\xi)=\hat\epsilon(w_i^p\xi)\pm\hat\epsilon(u_{i,1}d\xi)=0$.
\item
$\hat\epsilon(dv_{i,1})=\hat\epsilon((c_i-\gamma_i)w_i^{p-1})=0$ because $\gamma_i\in T$.
\end{itemize}

This proves that $\hat\epsilon(d\hat A)=0$, in other words
$\hat\epsilon$ is a chain map. That it induces an
epimorphism in cohomology in degree $n$ follows immediately from the
facts that $\epsilon$ does and that
$\epsilon=\hat\epsilon j$.
\end{proof}
This completes our construction of an oriented CDGA $(\hat A,\hat
d,\hat \epsilon)$ quasi-isomorphic to $(A,d,\epsilon)$.

\section{Extending the range of half-acyclicity}
The aim of this section is to prove that the construction of the
previous section increases the range in which the set of orphans is
half-acyclic. More precisely we will prove the following:
\begin{prop}
\label{prop-indstep} Let $(A,d,\epsilon)$ be an oriented CDGA
satisfying the assumptions \refequ{equ-2red} and let $k\geq \hn+1$.
Then the CDGA $(\hat A,\hat d,\hat\epsilon)$ constructed in
the previous section also satisfies the assumptions
\refequ{equ-2red}.

Moreover if the set $\Or$ of
orphans of $(A,d,\epsilon)$ is $(k-1)$-half-acyclic, then the set $\hat\Or$ of orphans of $(\hat
A,\hat d,\hat \epsilon)$ is $k$-half-acyclic.
\end{prop}
The proof of this proposition consists of a long series of lemmas.
Recall the spaces $\Gamma$ from equation (\ref{equ-Gamma})
and $Z$ from above Lemma \ref{lemma-ZG}.

Notice that by assumption \ref{equ-2red}(iii), $n\geq 7$
and hence $k\geq 5$.

\begin{lemma}
\label{lemma-Oi} If $i>n-k+2$ then $\Or^i\subset \hat \Or^i$.
\end{lemma}
\begin{proof}
Since $n-i<k-2$ we have that $\hat A^{n-i}=A^{n-i}$. Therefore
$\hat\epsilon(\Or^i.\hat A^{n-i})=\epsilon(\Or^i. A^{n-i})=0$. So
$\Or^i\subset \hat \Or^i$.
\end{proof}
\begin{lemma}
\label{lemma-Oikerd} For $i=k-2$, $k-1$ or $k$, we have $\hat
\Or^i\cap\ker d\subset \Or^i \cap\ker d$.
\end{lemma}
\begin{proof}
\underline{Case 1: $i=k-2$.}

\noindent
We have $\hat
A^{k-2}=A^{k-2}\oplus\langle\{w_i\}\rangle$. Let
$\omega=\xi+\sum_i r_iw_i\in\hat A^{k-2}$ with $\xi\in A^{k-2}$
and $r_i\in\Bk$. Then
$$d\omega=(d\xi-\sum r_i\gamma_i)+\sum r_ic_i\in A^{k-1}\oplus
\langle\{c_i\}\rangle.
$$
Therefore if $d\omega=0$ then we must also have $r_i=0$ for
each $i$. This implies that $\hat A^{k-2}\cap\ker d\subset
A^{k-2}\cap\ker d$. Thus $\hat \Or^{k-2}\cap\ker
d\subset\Or^{k-2}\cap\ker d$.

\underline{Case 2: $i=k-1$.}

\noindent Since $A^1=0$ and $A^0=\Bk$,
$\hat A^{k-1}=A^{k-1}\oplus\langle\{c_i\}\rangle$.
Let $\omega=\xi+\sum
r_ic_i\in \hat A^{k-1}$ with
$\xi\in A^{k-1}$. Suppose that $\omega\in\hat\Or^{k-1}\cap\ker d$.
Then $\xi\in \Or^{k-1}$ because otherwise there would exist
$\xi^*\in A$ such that $\epsilon(\xi\xi^*)\not=0$, and since
$\hat\epsilon(c_i.A)=0$ we would have
$\hat\epsilon(\omega.\xi^*)\not=0$.

Also $d\omega=0$ implies that $\sum r_i\alpha_i=d(-\xi)\in
d(\Or^{k-1})$. But by definition of $\{\alpha_i\}$ we have
$\langle\{\alpha_i\}\rangle\cap d(\Or^{k-1})=0$. Therefore $r_i=0$
for each $i$, hence $\omega=\xi\in \Or^{k-1}\cap\ker d$.

\underline{Case 3: $i=k$.}

\noindent
Let $\{\lambda_j\}$ be a basis of $A^2$. By assumption
\ref{equ-2red}(ii) this basis consists of cocycles. Since
$n\geq 7$, we have $k-2>2$ and $\hat
A^k=A^k\oplus\langle\{w_i.\lambda_j\}\rangle$. Let
$\omega=\xi+\sum r_{ij}w_i\lambda_j\in\hat A^k$ with $\xi\in A^k$.
Then
$$d\omega=(d\xi-\sum r_{ij}\gamma_i\lambda_j)+\sum
r_{ij}c_i\lambda_j\in
A^{k+1}\oplus\langle\{c_i\cdot \lambda_j\}\rangle.$$ Therefore
$d\omega\not=0$ unless $r_{ij}=0$ for all $i,j$. This implies that
$\hat A^k\cap\ker d\subset A^k\cap\ker d$, hence $\hat
\Or^k\cap\ker d\subset \Or^k\cap\ker d$
\end{proof}
Now the rest of the proof of Proposition \ref{prop-indstep} splits
into three cases: $k=\hn+1$ and $n$ even, $k=(n+1)/2+1$ and $n$ odd,
and $k\geq\hn+2$.
\subsection{The case $n$ even and $k=\hn+1$.}
\begin{lemma}\label{perp}
Let
$$
\xymatrix
{
0 \ar[r] & A \ar[r]^i & B \ar[r]_p & C \ar[r]\ar@/_/[l]_r & 0
}
$$
be a short exact sequence of vector spaces,
$r\colon C\rightarrow B$ be a linear spliting of
$p$ and
$\langle \_, \_ \rangle\colon B\otimes B \rightarrow \Bk$ be a non-degenerate
bilinear form on $B$. If $\langle \im r, \im i \rangle=0$,
then $\langle r\_, r\_ \rangle\colon C\otimes C \rightarrow \Bk$
is a non-degenerate bilinear form on $C$.
\end{lemma}
\begin{proof}
For any $\gamma \in C\setminus \{ 0\}$,
there is a $b\in B$ such that $\langle r\gamma, b\rangle\not=0$. Thus
$\langle r\gamma, rpb\rangle\not=0$, since $(rpb)-b\in \im i$.
\end{proof}

Recall the space $\Gamma=\langle \{ \gamma_i \} \rangle$ defined in \refequ{equ-Gamma}.
\begin{lemma}
\label{lemma-GPD} If $n$ is even and $k=\hn+1$ then  the bilinear
form
$$
\Gamma\otimes\Gamma\to\Bk\,,\,\gamma\otimes\gamma'\mapsto\epsilon(\gamma.\gamma')
$$
is non degenerate.
\end{lemma}
\begin{proof}
Set $n=2m$ and $k=m+1$. As in the proof of Proposition
\ref{prop-piqi}, the short exact sequence $0\to\Or\to A\to
\bar A:=A/\Or\to0$ induces a short exact sequence
$$\xymatrix{
0\ar[r]& \Ho^m(A)\ar[r]^{\pi^*}&\Ho^m(\bar A)\ar[r]^-{\delta} &
\Ho^{m+1}(\Or)\ar[r]&0 }
$$
where $\delta$ is the connecting homomorphism.
Since $\Or^{m+1}\cap\ker
d=d(\Or^m)\oplus\langle\{\alpha_i\}\rangle$, we get that
$\Ho^{m+1}(\Or)=\langle\{[\alpha_i]\}\rangle$.
Let $[\bar\gamma_i]\in \Ho^m(\bar A)$ be the cohomology classes represented by
$\bar\gamma_i=\gamma_i\mod\Or\in A^m/\Or$.

By Proposition \ref{prop-nada}, $\epsilon$
induces a non-degenerate pairing $\langle \cdot, \cdot \rangle$
on $\Ho^n(\overline A)$. Let $[\alpha_i] \mapsto [\bar\gamma_i]$ define a linear section $r$
of $\delta$. By Lemma \ref{lemma-G.kerd} we have
$\epsilon(\Gamma\cdot \ker d)=0$, and hence $\langle \im r, \im \pi^* \rangle=0$. Thus
by Lemma \ref{perp} the pairing restricts to a non-degenerate pairing on $\im r$. Finally observe that
under the identification of $\im r$ with $\Gamma$ which sends $[\bar\gamma_i]$ to $\gamma_i$
the restricted pairing is sent to the pairing given in the statement of the lemma.
\end{proof}
\begin{lemma}
\label{lemma-alphanotOr} If $n$ is even and $k=\hn+1$ then
$\langle\{\alpha_i\}\rangle\cap\hat \Or^k=0$.
\end{lemma}
\begin{proof}
Let $\alpha:=\sum r_i\alpha_i\in\langle\{\alpha_i\}\rangle$. If
the $r_i$ are not all zero then by Lemma \ref{lemma-GPD} there
exist $r^*_j\in\Bk$ such that $\epsilon\left((\sum
r_i\gamma_i)(\sum r^*_j\gamma_j)\right)\not=0$. Then
$$\hat\epsilon\left((\sum
r_i\alpha_i)(\sum r^*_jw_j)\right)= \sum
r_ir^*_j\hat\epsilon((d\gamma_i).w_j)= \pm \sum
r_ir^*_j\hat\epsilon(\gamma_i.\gamma_j)\not=0.$$ Hence
$\alpha\not\in\hat\Or$ if $\alpha\not=0$.
\end{proof}

\begin{lemma}
\label{lemma-Ok-1even} If $n$ is even and $k=\hn+1$ then
$\Or^{k-1}\subset\hat\Or^{k-1}$
\end{lemma}
\begin{proof}
Let $\beta\in\Or^{k-1}$. Then $\hat
A^{n-(k-1)}=A^{k-1}\oplus \langle\{c_i\}\rangle$. Let
$\omega=\xi+\sum r_ic_i$ with $\xi\in A^{k-1}$. Then
$\hat\epsilon(\beta\omega)=\epsilon(\beta\xi)+\sum
r_i\hat\epsilon(\beta c_i)=0$. Therefore $\beta\in\hat\Or^{k-1}$.
So $\Or^{k-1}\subset\hat\Or^{k-1}$ and we are done.
\end{proof}

\begin{lemma}
\label{lemma-indneven} If $n$ is even and  $k=\hn+1$ then
$\hat\Or$ is $k$-half-acyclic.
\end{lemma}
\begin{proof}
We only need to check that $\hat\Or^k\cap \ker d\subset d(\hat\Or^{k-1})$.
By Lemma \ref{lemma-Oikerd}
$$\hat \Or^k\cap\ker d\subset \Or^k\cap\ker
d=d(\Or^{k-1})\oplus\langle\{\alpha_i\}\rangle.$$ By Lemma
\ref{lemma-Ok-1even} this implies that
\begin{equation}\label{equ-indneven}
\hat\Or^k\cap\ker d \subset
d(\hat\Or^{k-1})\oplus\langle\{\alpha_i\}\rangle.
\end{equation}
Since the set of orphans is a differential ideal, we also have $d(\hat\Or^{k-1})\subset\hat\Or^k\cap\ker
d$. This combined with Lemma \ref{lemma-alphanotOr} and inclusion
\refequ{equ-indneven} implies that $\hat\Or^k\cap \ker d\subset
d(\hat\Or^{k-1})$.
\end{proof}
\subsection{The case $n$ odd and $k=(n+1)/2+1$.}
Recall the space $Z$ defined before Lemma \ref{lemma-ZG}
\begin{lemma}
\label{lemma-Zci} If $n$ is odd and $k=(n+1)/2+1$ then
$Z^{k-1}\oplus\langle\{c_i\}\rangle\subset\hat\Or^{k-1}$.
\end{lemma}
\begin{proof}
Notice that $n-(k-1)=k-2$ and $\hat
A^{k-2}=A^{k-2}\oplus\langle\{w_i\}\rangle$. It is immediate to
check, using the definition \refequ{equ-epshat} of $\hat\epsilon$ and
the fact that $Z\subset T\cap \Or$,
that $\hat\epsilon(\hat A^{k-2}\cdot Z^{k-1})=0$. Also $\hat\epsilon(\hat
A^{k-2}\cdot c_i)=0$ since $c_j\not\in \hat A^{k-2}$.
\end{proof}
\begin{lemma}
\label{lemma-indnodd} If $n$ is odd and $k=(n+1)/2+1$ then $\hat\Or$
is $k$-half-acyclic.
\end{lemma}
\begin{proof}
We only need to check that $\hat\Or^k\cap\ker d\subset
d(\hat\Or^{k-1})$. Using Lemmas \ref{lemma-Oikerd} and \ref{lemma-ZG}
we have that
$$
\hat\Or^k\cap\ker d\subset\Or^k\cap\ker
d=d(\Or^{k-1})\oplus\langle\{\alpha_i\}\rangle= d(Z^{k-1})\oplus
d(\langle\{c_i\}\rangle)$$
By Lemma \ref{lemma-Zci} the
last set is included in  $d(\hat\Or^{k-1})$.
\end{proof}
\subsection{The case $k\geq\hn+2$.}
\begin{lemma}
\label{lemma-Oi=}
 If $\hn\leq i\leq k-3$ then $\hat \Or^i= \Or^i$.
\end{lemma}
\begin{proof}
If $i\leq k-3$ then $\hat A^i=A^i$, so $\hat \Or^i\subset \Or^i$.

If $\hn\leq i\leq k-3$ then
$$
i\geq\hn= n-\hn\geq n-(k-3)>n-k+2
$$
and by Lemma \ref{lemma-Oi} $\Or^i\subset\hat\Or^i$.
\end{proof}

\begin{lemma}
\label{lemma-Zk-2} If $k\geq \hn+2$ then
$Z^{k-2}\subset\hat\Or^{k-2}$.
\end{lemma}
\begin{proof}
First suppose that $n$ is odd or that $k\geq\hn+3$. In these cases
$2k>n+4$, hence $k-2>n-k+2$ which implies by Lemma \ref{lemma-Oi}
that $\Or^{k-2}\subset\hat\Or^{k-2}$. Since
$Z^{k-2}\subset\Or^{k-2}$, this completes the proof of the lemma in these cases.

Now suppose that $n$ is even and $k=\hn+2$, then $n-(k-2)=k-2$.
Since $Z\subset \Or$ we have $\epsilon(Z^{k-2}\cdot A^{k-2})=0$. Also
by definition of $\hat\epsilon$ since $Z\subset T$,
$\hat\epsilon(Z^{k-2}\cdot w_i)=0$. Since $\hat
A^{n-(k-2)}=A^{n-(k-2)}\oplus\langle\{w_i\}\rangle$ this implies
that $Z^{k-2}\subset\hat\Or^{k-2}$.
\end{proof}

\begin{lemma}
\label{lemma-k-1half} If $k\geq\hn+2$ and $\Or$ is
$(k-1)$-half-acyclic then so is $\hat\Or$.
\end{lemma}
\begin{proof}
For $\hn+1\leq i\leq k-3$, using Lemma \ref{lemma-Oi=} twice, we
get that
$$\hat\Or^i\cap\ker d= \Or^i\cap\ker d\subset d(\Or ^{i-1})=d(\hat\Or^{i-1}).
$$

By Lemmas \ref{lemma-Oikerd},
\ref{lemma-ZG}  and \ref{lemma-Zk-2} we have
$$\hat\Or^{k-1}\cap\ker d
\subset\Or^{k-1}\cap\ker d\subset d(\Or^{k-2})=d(Z^{k-2})\subset
d(\hat\Or^{k-2}).$$

Suppose that $k\geq\hn+3$ (otherwise there is no need to check
$(k-2)$-half-acyclicity.) By Lemma \ref{lemma-Oi=} $\hat
\Or^{k-3}= \Or^{k-3}$ and so by Lemma \ref{lemma-Oikerd} we have
$$\hat\Or^{k-2}\cap\ker d\subset\Or^{k-2}\cap\ker d\subset
d(\Or^{k-3})=d(\hat\Or^{k-3}).$$
\end{proof}

\begin{lemma}
\label{lemma-ciOk-1} If $k\geq\hn+2$ then
$\langle\{c_i\}\rangle\subset\hat\Or^{k-1}$.
\end{lemma}
\begin{proof}
By definition of $\hat\epsilon$ the only products with $c_i$ which
could prevent them from being orphans are
$$
\hat\epsilon(c_ic_j)=-\epsilon( \gamma_i\gamma_j)
$$
but those are zeros for degree reasons.
\end{proof}

\begin{lemma}
\label{lemma-indkhn+2}
If $k\geq\hn+2$ and $\Or$ is
$(k-1)$-half-acyclic then
 $\hat\Or$ is $k$-half-acyclic.
\end{lemma}
\begin{proof}
We already know by Lemma \ref{lemma-k-1half} that $\hat\Or$ is
$(k-1)$-half-acyclic. Since $k\geq\hn+2$ we have $k-1>n-k+2$ and
Lemma \ref{lemma-Oi} implies that $\Or^{k-1}\subset\hat\Or^{k-1}$.

By Lemma \ref{lemma-ciOk-1} and the definitions of $dc_i$, we have
$$\langle\{\alpha_i\}\rangle= d(\langle\{c_i\}\rangle)\subset
d(\hat\Or^{k-1}).$$
 Using Lemma \ref{lemma-Oikerd} and the
definition of $\{\alpha_i\}$ we get
$$
\hat\Or^k\cap\ker d\subset \Or^k\cap\ker d=d(\Or^{k-1})\oplus
\langle\{\alpha_i\}\rangle\subset d(\hat\Or^{k-1}).$$ This
proves that $\hat\Or$ is $k$-half-acyclic.
\end{proof}
\subsection{End of the proof of Proposition \ref{prop-indstep}}
\begin{proof}[Proof of Proposition \ref{prop-indstep}]
Since $n\geq 7$ we have $k\geq 5$ and also using the fact that
$j\colon A\to\hat A$ is a quasi-isomorphism by Lemma
\ref{lemma-jqi}, and that $\hat \epsilon$ is an orientation by
Lemma \ref{lemma-ori}, it is immediate to check that $(\hat A,\hat
d,\hat \epsilon)$ satisfies the assumptions \refequ{equ-2red}.

If $\Or$ is $(k-1)$-half-acyclic for some $k\geq\hn+1$ then
Lemmas \ref{lemma-indneven}, \ref{lemma-indnodd}, and
\ref{lemma-indkhn+2} imply that $\hat\Or$ is $k$-half-acyclic.
\end{proof}
\section{Proof of Theorem \ref{thm-main}}
We conclude the proof of our main theorem.
\begin{proof}[Proof of Theorem \ref{thm-main}]
If $n\leq 6$ then by \cite{Neisendorfer} the CDGA $(A,d)$ is
formal and we can take its cohomology algebra as the Poincar\'e
duality model.

Suppose that $n\geq 7$. Since $H(A,d)$ is simply-connected, by taking
a minimal Sullivan model we can suppose that $A$ is of finite
type,
 $A^0=\Bk$, $A^1=0$, and $A^2\subset\ker d$. Also there exists
 a chain map $\epsilon\colon A\rightarrow s^{-n} \Bk$ inducing a surjection in
 homology. So all the assumptions
\refequ{equ-2red} are satisfied. Taking $k=\hn+1$ if $n$ is even
or $k=(n+1)/2+1$ if $n$ is odd, the set of
orphans $\Or$ is $(k-1)$-half-acyclic because this condition is
empty. An obvious induction using Proposition \ref{prop-indstep}
yields a quasi-isomorphic oriented  model $\hat A$ for which
the set of orphans $\hat \Or$ is $(n+1)$-half-acyclic. Propositions
\ref{prop-nada} and
\ref{prop-piqi} imply that the quotient $A'=\hat A/\hat\Or$ is a
Poincar\'e CDGA quasi-isomorphic to $A$.
\end{proof}

\section{Equivalence of differential Poincar\'e duality algebras}
\label{sec-equ}

The next theorem shows that if we have two $3$-connected quasi-isomorphic
differential Poincar\'e duality algebras then they can be connected by quasi-isomorphisms involving
only differential Poincar\'e duality algebras.

\begin{thm}\label{thm-qi}
Suppose $A$ and $B$ are quasi-isomorphic finite type differential
Poincar\'e duality algebras of dimension at
least $7$ such that
$\Ho^{\leq 3}(A)=\Ho^{\leq 3}(B)=\Bk$ and $A^{\leq 2}=B^{\leq 2}=\Bk$.
There exists a differential Poincar\'e duality algebra $C$ and quasi-isomorphisms
$A \rightarrow C$ and $B\rightarrow C$.
\end{thm}
\begin{proof}
Let $\wedge V$ be a minimal Sullivan model of $A$. Then there exist quasi-isomorphisms
$f\colon \wedge V\rightarrow A$ and $g\colon \wedge V\rightarrow B$. Consider the factorization of
the multiplication map $\phi\colon\wedge V \otimes \wedge V\rightarrow \wedge V$ into a cofibration
$i\colon \wedge V\otimes \wedge V\rightarrow \wedge V\otimes \wedge V\otimes \wedge U$ followed by
a quasi-isomorphism $p\colon\wedge V\otimes \wedge V\otimes \wedge U \rightarrow \wedge V$.
Since $\Ho^{\leq 3}(\wedge V)=\Ho^{\leq 3}(A)=\Bk$, we can assume that $U^{\leq 2}=0$ and that
$U$ is of finite type.
Next consider the following diagram in which $C'$ is defined to make the bottom square a pushout.
$$
\xymatrix
{
& \wedge V \ar[r]^f \ar[d]_{in_1} & A \ar[d]^{in_1}\\
& \wedge V\otimes \wedge V \ar[dl]_{\phi}
\ar[r]^-{f\otimes g} \ar[d]_i & A\otimes B \ar[d]^h \\
\wedge V  & \wedge V\otimes \wedge V\otimes \wedge U\ar[l]^-p  \ar[r]_-k & C'
}
$$
The maps $in_1$ denote inclusion into the first factor. Since $p$ is a quasi-isomorphism and
$\phi\circ in_1=id$, $i\circ in_1$ is a quasi-isomorphism. Also since
$f\otimes g$ is a quasi-isomorphism, $i$ is a cofibration, the bottom square is a pushout,
and the properness of CDGA \cite[Lemma 8.13]{Baues}, $k$ is a quasi-isomorphism.
Finally since $f$, $k$ and
$i\circ in_1$ are quasi-isomorphisms, so too must $h\circ in_1$  be a quasi-isomorphism.

Since $A^{\leq 2}=B^{\leq 2}=\Bk$,
the algebra $A\otimes B$ satisfies (i) and (ii) of (\ref{equ-2red}).
Also $U$ is of finite type with $U^{\leq 2}=0$, so $C'=A\otimes B \otimes \wedge U$ satisfies
(i) and (ii) of (\ref{equ-2red}). Since $C'$ is
quasi-isomorphic to $A$, it also satisfies (iii), and
we can let $\epsilon\colon C' \rightarrow s^n\Bk$ be any orientation.
Next by using induction and Propositions
\ref{prop-indstep}, \ref{prop-nada} and \ref{prop-piqi}
we get a quasi-isomorphism $l\colon C'\rightarrow C$ such that $C$ is a
differential Poincar\'e duality algebra. Clearly the map $l\circ h\circ in_1\colon A\rightarrow C$ is a
quasi-isomorphism, and similarly $l\circ h\circ in_2\colon B\rightarrow C$ is a quasi-isomorphism,
thus completing the proof of the theorem.
\end{proof}

\noindent {\bf Conjecture:}
The hypotheses $\Ho^{\leq 3}(A)=\Bk$, $A^{\leq 2}=B^{\leq 2}=\Bk$ and dimension of $A$
at least $7$ in Theorem \ref{thm-qi} can be removed.
%


\begin{thebibliography}{xx}
\bibitem{ALH}
{\bf M. Aubry, J.-M. Lemaire, and S. Halperin} {\em Poincar\'e duality models},
preprint.
%
\bibitem{Baues}
{\bf H.J. Baues}
{\em Algebraic Homotopy}, Cambridge University Press (1989).
%

\bibitem{ChasSullivan}
{\bf M. Chas and D. Sullivan}
{\em String topology}
{preprint arXiv: GT/991159}

\bibitem{FHT}
{\bf Y. F\'elix, S. Halperin and J.-C. Thomas}
{\em Rational homotopy theory}, Graduate Texts in Mathematics, vol. 210, Springer-Verlig (2001).
%
\bibitem{FelixThomas}
{\bf Y. F\'elix and J.-C. Thomas}
{\em Rational BV-algebra in string topology.} Preprint arXiv:0507.4194, (2007).
%
%
\bibitem{FultonMacPherson}
{\bf W. Fulton and R. Mac Pherson} {\em A compactification of
configuration spaces.} Annals of Math. {\bf 139} (1994), 183--225.
%
\bibitem{Kriz}
{\bf I. Kriz} {\em On the rational homotopy type of configuration spaces.} Annals of Math.
{\bf 139} (1994), 227--237.
%
\bibitem{Lambrechts-cochainthick}
{\bf P. Lambrechts} {\em Cochain model for thickenings and its
application to rational LS-category.} Manuscripta Math. 103
(2000), 143--160.
%
\bibitem{LS-FM2}
{\bf P. Lambrechts and D. Stanley}
{\em The rational homotopy type of configuration spaces of two points.}
Annales Inst. Fourier {\bf 54} (2004), 1029--1052.
%
\bibitem{LS-dgmodFMk}
{\bf P. Lambrechts and D. Stanley} {\em A remarkable DG-module model for configuration spaces.}
 Mittag-leffler preprint (spring 2006, vol.38). Submitted. Preprint arXiv:math.0707.2350.
%
%
\bibitem{Menichi}
{\bf L. Menichi} {\em Batalin-Vilkovisky algebra
structures on Hochschild cohomology.} Preprint arXiv:math.QA/0711.1946, (2007).
%
\bibitem{Neisendorfer}
{\bf J. Neisendorfer and T. Miller} {\em Formal and coformal spaces.}
Illinois J. Math. {\bf 22} (1978), 565--580.
%
\bibitem{Stasheff}
{\bf J. Stasheff}
{\em Rational Poincar\'e duality spaces. }
Illinois J. Math. {\bf 27} (1983), 104--109.
%
\bibitem{Tradler-BV}
{\bf T. Tradler}
{\em The BV Algebra on Hochschild Cohomology Induced by Infinity
        Inner Products},
Preprint arXiv:math.QA/0210150.
%
\bibitem{Sul:inf}
{\bf D. Sullivan},
{\em Infinitesimal computations in topology},
Inst. Hautes \'Etudes Sci. Publ. Math. No. {\bf 47}, (1977), 269--331.
%
\bibitem{Yang}
{\bf T. Yang}, {\em A Batalin-Vilkovisky algebra structure on the Hochschild cohomology of
truncated polynomials.}
Master's Thesis, University of Regina, (2007).

\end{thebibliography}
\end{document}